\documentclass[final,3p,11pt,sort&compress]{elsarticle}
\usepackage{amsmath}
\usepackage{amsfonts}
\usepackage{amssymb}
\usepackage{graphicx}
\usepackage{subfigure}
\usepackage{threeparttable}
\usepackage{xspace}

\newcommand{\ignore}[1]{}
\newcommand{\etal}[0]{{\em et~al.\@}\xspace}
\newcommand{\eg}[0]{{e.g.\@}\xspace}
\newcommand{\ie}[0]{{i.e.\@}\xspace}

\newcommand{\fnc}[1]{\ensuremath{\mathcal{#1}}}
\newcommand{\bsym}[1]{\ensuremath{\boldsymbol{#1}}}

\DeclareMathOperator{\diag}{diag}

\newtheorem{thm}{Theorem}
\newtheorem{lem}{Lemma}
\newtheorem{prop}{Proposition}
\newtheorem{coro}{Corollary}
\newdefinition{define}{Definition}
\newproof{pf}{Proof}

\begin{document}

\begin{frontmatter}

\title{Summation-By-Parts Operators and High-Order Quadrature\tnoteref{t1}}

\tnotetext[t1]{This work was supported by the Natural Sciences and
Engineering Research Council (NSERC), the Canada Research Chairs
program, Bombardier Aerospace, Mathematics of Information Technology
and Complex Systems (MITACS), and the University of Toronto}

\author[utias]{J.E.~Hicken\corref{cor1}\fnref{fn1}}
\ead{jehicken@stanford.edu}
\author[utias]{D.W.~Zingg\fnref{fn2}}
\ead{dwz@oddjob.utias.utoronto.ca}

\cortext[cor1]{corresponding author}

\fntext[fn1]{Postdoctoral fellow, Stanford University} 
\fntext[fn2]{Professor and Director, Tier 1 Canada Research Chair
in Computational Aerodynamics, J.\ Armand Bombardier Foundation Chair
in Aerospace Flight}

\address[utias]{Institute for Aerospace Studies, University of
  Toronto, Toronto, Ontario, M3H 5T6, Canada}

\begin{abstract}
Summation-by-parts (SBP) operators are finite-difference operators
that mimic integration by parts.  This property can be useful in
constructing energy-stable discretizations of partial differential
equations.  SBP operators are defined by a weight matrix and a
difference operator, with the latter designed to approximate $d/dx$ to
a specified order of accuracy.  The accuracy of the weight matrix as a
quadrature rule is not explicitly part of the SBP definition.  We show
that SBP weight matrices are related to trapezoid rules with end
corrections whose accuracy matches the corresponding difference
operator at internal nodes.  The accuracy of SBP quadrature extends to
curvilinear domains provided the Jacobian is approximated with the
same SBP operator used for the quadrature.  This quadrature has
significant implications for SBP-based discretizations; for example,
the discrete norm accurately approximates the $L^{2}$ norm for
functions, and multi-dimensional SBP discretizations accurately mimic
the divergence theorem.
\end{abstract}


\end{frontmatter}

\section{Introduction}

Partial differential equations (PDEs) are often solved numerically in
order to approximate a functional that depends on the solution; for
example, when computational fluid dynamics is used to estimate the
lift and drag on an aerodynamic body.  For integral functionals, such
as lift and drag, a quadrature rule may be needed to numerically
integrate the discrete solution.  When we are free to choose the
quadrature weights and abscissas, Guassian quadrature is often the
optimal choice.  However, the choice of quadrature rule is less clear
for the uniform grids that arise in finite-difference methods.

In this paper, we investigate a quadrature rule that is particularly
well suited for high-order summation-by-parts (SBP) finite-difference
methods \cite{kreiss:1974}.  SBP operators lead to linearly
time-stable discretizations of well-posed PDEs, and they have been
used to construct efficient discretizations of the Euler
\cite{mattsson:2004,hicken:2008}, Navier-Stokes
\cite{mattsson:2008,nordstrom:2009,osusky:asm2010}, and Einstein
equations \cite{pazos:2009}.  The high-order quadrature in question is
based on the weight matrix that forms part of the definition of SBP
operators.  This result is somewhat surprising, because the accuracy
of the quadrature induced by the weight matrix is not explicitly part
of the SBP definition.  To our knowledge, the relationship between SBP
operators and quadrature has not been discussed previously in the
literature.

In the context of high-order finite-difference methods, including
those based on SBP operators, several classical quadrature rules are
available to accurately evaluate integral functionals; for example,
composite Newton-Cotes rules and Gregory-type formulae.  Why introduce
a new quadrature rule based on SBP weight matrices?  While accuracy is
important, we may also want the functional estimate to obey some
property or properties of the true functional, and this is the value
of SBP-based quadrature.

Consider the volume integral of the divergence of a vector field over
a compact domain.  The resulting functional is equivalent to the flux
of the vector field over the domain's boundary, in light of the
divergence theorem.  This is a fundamental property of the functional
that we may want a discretization and quadrature to preserve.  We say
a functional estimate respects, or mimics, the divergence theorem if
1) it is accurate, and 2) the discrete quadrature over the volume
produces a discrete quadrature over the surface.

In general, classical quadrature rules for uniformly spaced data will
not mimic the divergence theorem in the above sense when applied to an
arbitrary high-order finite-difference approximation of the
divergence; typically, they will satisfy the first but not the second
property.  In contrast, we will show that an SBP discretization does
mimic the divergence theorem when numerically integrated using its
corresponding weight matrix.

\ignore{
In contrast, summation-by-parts (SBP) operators are finite-difference
approximations that are designed explicitly to satisfy the second
property, but not necessarily the first \cite{kreiss:1974}.  One can
show that the second property leads to linearly time-stable
discretizations of well-posed PDEs.  SBP operators have been used to
construct efficient discretizations of the Euler
\cite{mattsson:2004,hicken:2008}, Navier-Stokes
\cite{mattsson:2008,nordstrom:2009,osusky:asm2010}, and Einstein
equations \cite{pazos:2009}.

SBP operators are defined by a difference operator and a weight
matrix.  The weight matrix is used to construct a discrete norm for
stability analysis, and it provides the necessary flexibility to meet
both property 2 and the accuracy requirements of the SBP
finite-difference operator.

In this paper, we show that SBP weight matrices define quadrature
rules with the same order of accuracy as the corresponding
finite-difference operator at internal nodes.  Consequently, the
weight matrices provide a natural quadrature rule for SBP
discretizations and related functionals.  In particular, an SBP
discretization will mimic the divergence theorem: they obey both
properties 1 and 2 above.  This result is somewhat surprising,
because the accuracy of the quadrature induced by the weight matrix is
not explicitly part of the SBP definition.  To our knowledge, the
relationship between SBP operators and quadrature has not been
discussed in the literature.
}

The paper is organized as follows.  Section \ref{sec:define}
introduces notation and formally defines SBP operators.
Section~\ref{sec:theorem} presents the main theoretical results.  In
particular, we derive conditions on the quadrature weights for the
class of trapezoid rules with end corrections.  These conditions are
used to establish the accuracy of SBP-based quadrature.  Subsequently,
we consider the impact of coordinate transformations on SBP quadrature
and show that the quadrature remains accurate on curvilinear
multi-dimensional domains. In Section~\ref{sec:example} we verify the
theoretical results with several numerical examples.  The implications
of SBP quadrature are summarized and discussed in
Section~\ref{sec:conclude}.

\section{Notation and definitions}\label{sec:define}

We try to remain consistent with the notation used by Kreiss and
Scherer in their original work \cite{kreiss:1974}, as well as Strand's
subsequent work \cite{strand:1994}.

The interval $[0,1]$ is partitioned into $n+1$ evenly spaced points
$x_{v} = vh, \; v = 0,1,\ldots,n$, with mesh spacing $h = 1/n$.
Finite intervals other than $[0,1]$, as well as nonuniform node
spacing, can be accommodated by introducing an appropriate mapping
(see Section \ref{sec:mapping}).  For arbitrary $\fnc{U}(x) \in
C^{p}[0,1]$, we use $u_{v} = \fnc{U}(x_{v})$ to denote the restriction
of $\fnc{U}$ to the grid $x_{v}$.

\begin{define}[Summation-By-Parts Operator]\label{def:sbp}
  The matrix $D \in \mathbb{R}^{(n+1)\times (n+1)}$ is a
  summation-by-parts operator for the first derivative on the mesh
  $\{x_{v}\}_{v=0}^{n}$ if it has the form
  \begin{equation*}
    D = H^{-1} Q,
  \end{equation*}
  where the weight matrix $H \in \mathbb{R}^{(n+1)\times (n+1)}$ is a
  symmetric-positive-definite matrix, and $Q \in
  \mathbb{R}^{(n+1)\times (n+1)}$ satisfies
  \begin{equation*}
    Q + Q^{T} = \diag{(-1,0,0,\ldots,1)}.
  \end{equation*}
  Furthermore, the truncation error of the difference operator $D$ in
  approximating $d/dx$ is order $h^{2s}$ at the internal nodes,
  $\{x_{v}\}_{v=r}^{n-r}$, and order $h^{\tau}$ at the boundary nodes,
  $\{x_{v}\}_{v=0}^{r-1}$ and $\{x_{v}\}_{v=n-r+1}^{n}$, where $\tau, r,
  s \geq 1$.
\end{define}

In other words, the SBP operator $D$ approximates $d/dx$ and has a
particular structure.  In general, the order of accuracy of the
difference stencil at internal nodes is different than the order of
accuracy of the stencil at boundary nodes.  The even order of accuracy
$2s$ for the internal nodes is a consequence of using
centered-difference schemes, which provide the lowest error for a
given stencil size.  For a $2s$-order accurate scheme, the derivative
at the internal nodes is approximated as
\begin{align*}
  \frac{d \fnc{U}}{dx}(x_{w}) &\approx \sum_{v=1}^{s} \frac{\alpha_{v}}{h} 
  (u_{w+v} - u_{w-v}), \qquad r \leq w \leq n-r, \\
\intertext{where the coefficients are defined by 
  (see \cite{li:2005}, for example)}
  \alpha_{v} &= \frac{(-1)^{v+1} (s!)^{2}}{v (s+v)! (s-v)!}.
\end{align*}
The following lemma from \cite{strand:1994} lists some identities
that the $\alpha_{v}$ satisfy; these identities will be useful in our
subsequent analysis.

\begin{lem}\label{lem:alpha}
The coefficients $\alpha_{v}$ that define a $2s$-order accurate SBP
operator at internal nodes satisfy
\begin{align*}
  \sum_{v=1}^{s} \alpha_{v} v^{2j+1} = 
  \begin{cases} 
    \frac{1}{2}, & j = 0, \\
    0,           & j = 1,2,\ldots,s-1.
  \end{cases}   
\end{align*}
\end{lem}

We turn our attention to the weight matrix $H$, which is the focus of
this paper.  Since $H$ is symmetric-positive-definite, we can use it
to define an inner product and corresponding norm for vectors.  Let
$u, z \in \mathbb{R}^{n+1}$ be two discrete functions on the grid
nodes, \ie $u_{v} = \fnc{U}(x_{v})$ and $z_{v} = \fnc{Z}(x_{v})$.
Then
\begin{equation*}
  (u,z)_{H} \equiv u^{T} H z, \quad \text{and } \quad  
  \| u \|_{H}^{2} \equiv (u,u)_{H},
\end{equation*}
define the $H$ inner product and $H$ norm, respectively.  Using the
SBP-operator definition and the $H$ inner product, we have
\begin{equation}\label{eq:SBP}
  (u,D\,z)_{H} = -(D\,u,z)_{H} - u_{0} z_{0} + u_{n} z_{n}.
\end{equation}
Equation~\eqref{eq:SBP} expresses the fundamental property of SBP operators
and is the discrete analog of 
\begin{equation}\label{eq:IBP}
  \int_{a}^{b} \fnc{U} \frac{d \fnc{Z}}{dx} dx = -\int_{a}^{b} \fnc{Z} \frac{d \fnc{U}}{dx} dx + \left. \fnc{U}\fnc{Z} \right|_{x=a}^{x=b}.
\end{equation}
This property of SBP operators is what leads to energy-stable
discretizations of partial differential equations.  However, while
\eqref{eq:SBP} is analogous to integration by parts, it remains to be
shown that \eqref{eq:SBP} is an accurate discretization of
\eqref{eq:IBP}.

In this work, we will consider $H$ matrices with the block structure
\begin{equation}\label{eq:Hfull}
  H = h \begin{pmatrix} H_{L} & 0 & 0 \\ 0 & I & 0 \\ 0 & 0 & H_{R} \end{pmatrix},
\end{equation}
where $H_{L}, H_{R} \in \mathbb{R}^{r\times r}$ are
symmetric-positive-definite matrices.  Assuming that $H_{L}$ and
$H_{R}$ are dense matrices --- the so-called full-norm case --- Kreiss
and Scherer \cite{kreiss:1974} established the existence of SBP
operators that achieve an order of accuracy of $\tau = 2s-1$ at the
boundary with $r = 2s$.  Strand \cite{strand:1994} showed that $2s-1$
accuracy can be maintained at the boundary in the case of a
restricted-full norm, which uses
\begin{equation*}
  H_{L} = \begin{pmatrix} h_{00} & 0 \\ 0 & \bar{H}_{L} \end{pmatrix}
  \quad \text{and } \quad 
  H_{R} = \begin{pmatrix} \bar{H}_{R} & 0 \\ 0 & h_{00} \end{pmatrix}
\end{equation*}
with $\bar{H}_{L},\bar{H}_{R} \in \mathbb{R}^{(r-1)\times(r-1)}$ and
$r = 2s+1$.  

In general, SBP weight matrices of the form \eqref{eq:Hfull}
satisfy the compatibility conditions described in the following proposition
\cite{kreiss:1974}; these conditions will be used later to establish
the accuracy of quadrature rules based on full and restricted-full
$H$ matrices.

\begin{prop}\label{prop:fullH}
Let $H \in \mathbb{R}^{(n+1)\times (n+1)}$ be an SBP weight matrix
with the block structure \eqref{eq:Hfull}.  Then $H_{L}$ satisfies
\begin{equation*}
  j e_{i}^{T} H_{L} e_{j-1} + i e_{j}^{T} H_{L} e_{i-1}  
  = -(-r)^{i+j} + J_{i,j}, \qquad 0 \leq i,\,j \leq \tau,
\end{equation*}
where $e_{j}^{T} \equiv (-1)^{j}
\begin{pmatrix} r^{j} & (r-1)^{j} & \cdots & 1^{j} \end{pmatrix}$, 
with the convention $e_{-1} = 0$, and
\begin{equation*}
J_{i,j} = \sum_{v=1}^{s} \alpha_{v} \left[ \sum_{w=0}^{v-1} w^{j} (w -
v)^{i} + w^{i} (w - v)^{j} \right], \qquad i+j \geq 1.
\end{equation*}
\end{prop}

Kreiss and Scherer also showed that it is possible to define SBP
operators with diagonal $H$ matrices, \ie
\begin{align*}
H_{L} &= \diag{(\lambda_{0},\lambda_{1},\ldots,\lambda_{r-1})} \\
H_{R} &= \diag{(\lambda_{r-1},\ldots,\lambda_{1},\lambda_{0})}
\end{align*}
with $\lambda_{i} > 0$.  These ``diagonal norms'' are important
because, unlike full and restricted-full norms, they lead to provably
stable PDE discretizations on curvilinear grids \cite{svard:2004b}.
However, diagonal-norm SBP operators are limited to $\tau = s$
accuracy at the boundary when the internal accuracy is $2s$.
Consequently, the solution accuracy of hyperbolic systems discretized
with such SBP operators is limited to order $s+1$
\cite{gustafsson:1975}.  Nevertheless, one can show that functionals
based on the solution of dual consistent diagonal-norm SBP
discretizations are $2s$-order accurate \cite{hicken:cfdsc2010,
  hicken:supfun2011}.

When the weight matrix $H$ is diagonal, Kreiss and Scherer
\cite{kreiss:1974} showed that its elements are defined by the
relations in following proposition.
\begin{prop}\label{prop:diagH}
  Let $H \in \mathbb{R}^{(n+1)\times (n+1)}$ be a diagonal SBP weight
  matrix with $r = 2\tau = 2s$.  Then the diagonal elements $\lambda_{v}$ of
  $H_{L}$ and $H_{R}$ satisfy the relations
  \begin{equation*}
    j \sum_{v=0}^{r-1} \lambda_{v} (r - v)^{j-1} = \begin{cases} 
      \phantom{2 \sum_{v=1}^{s}} (r)^{j} - (-1)^{j} \beta_{j} ,   
      & j = 1,2,\ldots,2s-1 \\
      (r)^{2s} - 2 \sum_{v=1}^{s} \alpha_{v} 
      \sum_{w=0}^{v-1} w^{s}(w - v)^{s}, & j = 2s,
    \end{cases}
  \end{equation*}
  where $\beta_{j}$ is the $j^{\text{th}}$ Bernoulli number.
\end{prop}

\ignore{
Remark: would it be possible to increase the bandwidth of $H_{L}$ and $H_{R}$
such that the accuracy is retained?
}

\section{Theory}\label{sec:theorem}

\subsection{One-dimensional SBP Quadrature}

To establish the accuracy of SBP-based quadratures, we need the
following theorem that places constraints on the coefficients of a
certain class of quadrature rules for uniformly spaced data;
specifically, the trapezoid rule with end corrections.  The theorem is
a direct consequence of substituting finite-difference approximations
into the Euler-Maclaurin sum formula.

\begin{thm} \label{thm:eulermac} 
Consider a set of $n+1$ uniformly spaced points, $x_{v} = vh, \; v =
0,1,\ldots,n$, with constant mesh spacing $h = 1/n$.  A quadrature of the form
\begin{equation*}
  \fnc{I}(u) \equiv h \left( \sum_{v=0}^{r-1} \sigma_{v} u_{v} + \sum_{v=r}^{n-r} u_{v} + \sum_{v=0}^{r-1} \sigma_{v} u_{n-v} \right)
\end{equation*}
is a $q$-order accurate approximation of $\int_{0}^{1} \fnc{U}\, dx$
for $\fnc{U} \in C^{2m+2}[0,1]$, where $q-1 \leq r$ and $q \leq 2m + 2$,
if and only if the coefficients $\{ \sigma_{v} \}_{v=0}^{r-1}$ satisfy
\begin{equation}\label{eq:cnstrts}
  j \sum_{v=0}^{r-1} \sigma_{v} (r - v)^{j-1} = r^{j} - (-1)^{j} \beta_{j},
  \qquad j = 1,2,\ldots,q-1,
\end{equation}
\end{thm}
 
\begin{pf}
  Consider the Euler-Maclaurin sum formula applied to $\fnc{U}(x)$
  \cite{hildebrand:1974}:
  \begin{align}
    \int_{0}^{1} \fnc{U}(x) \, dx &= h \sum_{v=0}^{n} u_{v} 
    + \sum_{k=1}^{2m} \frac{\beta_{k}}{k!} h^{k} 
    \left( u^{(k-1)}_{0} - (-1)^{k}u^{(k-1)}_{n} \right) + E_{2m}, 
    \label{eq:eulermac} \\
    \intertext{where $u^{(k-1)}_{v} \equiv D^{(k-1)} \fnc{U}(x_{v})$,
      $2m < q \leq 2m+2$, and the error term is given by} E_{2m} &=
      \frac{\beta_{2m+2} h^{2m+2}}{(2m + 2)!} D^{(2m+2)} \fnc{U}(\xi),
      \notag
  \end{align}
  with $\xi \in (0,1)$.  Suppose the function derivatives at $x = 0$
  and $x = 1$ are replaced with finite-difference approximations
  involving the first $r$ and last $r$ internal points, respectively.
  Moreover, assume that the approximation to $u_{v}^{(k-1)}$ is accurate
  to $\text{O}(h^{q-k})$, where $q-1 \leq r$; consequently, the
  approximations are exact for polynomials up to at least degree
  $q-1$.  Let $\{ \delta_{v}^{(k-1)} \}_{v=0}^{r-1}$ denote the
  coefficients defining the finite-difference approximation of
  $u_{0}^{(k-1)}$, such that
  \begin{equation*}
    u^{(k-1)}_{0} = \sum_{v=0}^{r-1} \frac{\delta_{v}^{(k-1)}}{h^{k-1}} u_{v} 
    + \text{O}(h^{q - k}).
  \end{equation*}
  Substituting the finite-difference approximations into
  \eqref{eq:eulermac}, and noting that the coefficients for odd
  derivatives must be negated at $x=1$, we find
  \begin{align*}
    \int_{0}^{1} \fnc{U}(x) \, dx 
    &= h \sum_{v=0}^{n} u_{v} 
    + \sum_{k=1}^{2m} \frac{\beta_{k}}{k!} h^{k} 
      \sum_{v=0}^{r-1} \frac{\delta_{v}^{(k-1)}}{h^{k-1}} 
      \left( u_{v} + u_{n-v} \right) + \text{O}(h^{q}) + \text{O}(h^{2m+2}) \\
    &= h \left( \sum_{v=0}^{r-1} \sigma_{v} u_{v} 
    + \sum_{v=r}^{n-r} u_{v} + \sum_{v=0}^{r-1} \sigma_{v} u_{n-v} \right) 
    + \text{O}(h^{q})
  \end{align*}
  where
  \begin{equation}
    \sigma_{v} = 1 + \sum_{k=1}^{2m} \frac{\beta_{k}}{k!}
    \delta_{v}^{(k-1)}, \qquad v = 0,1,\ldots,r-1. \label{eq:parsol}
  \end{equation}
  Next, we will show that these $\sigma_{v}$ are the same ones that
  satisfy \eqref{eq:cnstrts}, a set of $q-1$ conditions that are
  independent of the $\delta_{v}^{(k)}$.  Substituting the above
  expression for $\sigma_{v}$ into \eqref{eq:cnstrts}, we find
  \begin{multline}\label{eq:cnstrts2}
    j \sum_{v=0}^{r-1} \sigma_{v} (r - v)^{j-1} 
    = j \sum_{v=0}^{r-1} (r-v)^{j-1} + j \sum_{k=1}^{2m} \frac{\beta_{k}}{k!} 
    \sum_{v=0}^{r-1} \delta_{v}^{(k-1)} (r - v)^{j-1}, \\
    j = 1,2,\ldots,q-1.
  \end{multline} 
  The first term on the right-hand side can be recast using the sum of
  powers formula\footnote{We use the sum of powers formula
  that is consistent with $\beta_{1} = -\frac{1}{2}$}:
  \begin{equation}\label{eq:term1}
    j \sum_{v=0}^{r-1} (r-v)^{j-1} = r^{j} + \sum_{k=1}^{j-1} (-1)^{k}
    \binom{j}{k} \beta_{k} r^{j-k}.
  \end{equation}
  For the second term, we recognize that $(r - v)^{j-1}$ is the
  discrete representation of the polynomial $p_{j-1}(x) \equiv h^{-(j-1)} (rh
  - x)^{j-1}$; therefore, since the finite-difference approximations
  are exact for polynomials of degree $q-1$, we have
  \begin{equation*}
    \sum_{v=0}^{r-1} \delta_{v}^{(k-1)} (r - v)^{j-1} = 
    \begin{cases}
      (-1)^{k-1} \frac{(j-1)!}{(j-k)!} r^{j-k}, & k \leq j \\
      0, & k > j,
    \end{cases}
  \end{equation*}
  and
  \begin{equation}\label{eq:term2}
    j \sum_{k=1}^{2m} \frac{\beta_{k}}{k!} 
    \sum_{v=0}^{r-1} \delta_{v}^{(k-1)} (r - v)^{j-1} 
    = \sum_{k=1}^{j} (-1)^{k-1} \binom{j}{k} \beta_{k} r^{j-k}.
  \end{equation}
  Substituting \eqref{eq:term1} and \eqref{eq:term2} into
  \eqref{eq:cnstrts2}, and recalling that the odd Bernoulli numbers
  greater than one are zero, we have
  \begin{align*}
    j \sum_{v=0}^{r-1} \sigma_{v} (r - v)^{j-1} 
    &= r^{j} + \sum_{k=1}^{j-1} (-1)^{k}\binom{j}{k} \beta_{k} r^{j-k} + 
    \sum_{k=1}^{j} (-1)^{k-1} \binom{j}{k} \beta_{k} r^{j-k} \\
    &= r^{j} - (-1)^{j} \beta_{j},
  \end{align*}
  for $j = 1,2,\ldots,q-1$.  Thus, we have shown that the $\sigma_{v}$
  satisfy \eqref{eq:cnstrts} when the quadrature is $q$-order accurate.

  We need the general solution of \eqref{eq:cnstrts} to show
  that these conditions are sufficient for the quadrature to be
  $q$-order accurate.  We have already shown that \eqref{eq:parsol} is a
  particular solution of the linear equations \eqref{eq:cnstrts}, so
  we need to determine the form of the homogeneous solution, \ie the
  null space of the matrix on the left side of \eqref{eq:cnstrts}.

  As noted above, $(r-v)^{j-1}$ is simply the polynomial $p_{j-1}(x) =
  h^{-(j-1)}(rh - x)^{j-1}$ evaluated at the nodes.  The derivative
  operator $D^{(k-1)}$ with $q \leq k \leq r$ will annihilate
  $p_{j-1}(x)$, since $j \leq q-1$; therefore, any finite difference
  approximation that is a consistent approximation of $h^{k-1}
  D^{(k-1)}$, $q \leq k \leq r$, will annihilate $p_{j-1}(x_{v}) =
  (r-v)^{j-1}$.  If we let $\{ \mu_{v}^{(k-1)} \}_{v=0}^{v=r-1}$
  denote the coefficients of such a finite difference approximation,
  then the general solution to \eqref{eq:cnstrts} can be written as
  \begin{equation}\label{eq:gensol}
    \sigma_{v} = 1 + \sum_{k=1}^{2m} \frac{\beta_{k}}{k!} \delta_{v}^{(k-1)}
    + \sum_{k=q}^{r} \gamma_{k-q} \mu_{v}^{(k-1)},    
  \end{equation}
  where $\{ \gamma_{0}, \gamma_{1}, \ldots, \gamma_{r-q} \}$
  parameterizes the null space.  When $r = q-1$, the null space is
  trivial, and the second sum does not appear in \eqref{eq:gensol}.

  Substituting the general solution into the quadrature yields
  \begin{align*}
    \fnc{I}(u) &= h \left( \sum_{v=0}^{r-1} \sigma_{v} u_{v} 
    + \sum_{v=r}^{n-r} u_{v} + \sum_{v=0}^{r-1} \sigma_{v} u_{n-v} \right) \\
    &= h \sum_{v=0}^{n} u_{v} 
    + \sum_{k=1}^{2m} \frac{\beta_{k}}{k!} h^{k} 
      \sum_{v=0}^{r-1} \frac{\delta_{v}^{(k-1)}}{h^{k-1}} 
      \left( u_{v} + u_{n-v} \right)  \\
    &\phantom{= h \sum_{v=0}^{n} u_{v}}\qquad
    + h \sum_{k=q}^{r} \gamma_{k-q} \sum_{v=0}^{r-1} 
    \mu_{v}^{(k-1)} (u_{v} + u_{n-v}) \\
    &= \int_{0}^{1} \fnc{U}(x) \, dx + \text{O}(h^{q})  + 
    \sum_{k=q}^{r} \gamma_{k-q}  h^{k} \left( u_{0}^{(k-1)} 
    + (-1)^{k} u_{n}^{(k-1)} \right)  \\
    &= \int_{0}^{1} \fnc{U}(x) \, dx + \text{O}(h^{q}).
  \end{align*}
  Therefore, we have shown that \eqref{eq:cnstrts} is sufficient for
  the quadrature to be $q$-order accurate, which completes the proof.
\end{pf}

If we choose $q-1 = r$, Theorem \ref{thm:eulermac} provides a closed
set of equations for constructing high-order quadrature rules for
uniformly spaced data with equal weights on the internal points.  More
generally, we may choose $q-1 < r$, in which case the additional
degrees of freedom can be used to achieve other objectives.  For
example, setting $\sigma_{0}$ to zero, so that only strictly internal
points are used.

Theorem \ref{thm:eulermac} encompasses many existing quadrature rules,
including the Gregory class of formulae, and it could be used to
construct an unlimited number of novel trapezoid rules with end
corrections.  However, our interest in Theorem \ref{thm:eulermac} is
not in constructing new quadrature rules, but in its consequences for
SBP weight matrices.

\begin{coro}\label{cor:sbpnorm}
Let $H$ be a full, restricted-full, or diagonal weight matrix from an
SBP first-derivative operator $D = (H^{-1}Q)$, which is a
$2s$-order-accurate approximation to $d/dx$ in the interior.  Then the
$H$ matrix constitutes a $2s$-order-accurate quadrature for integrands
$\fnc{U}\in C^{2s}$.
\end{coro} 
 
\begin{pf}
For diagonal SBP weight matrices the result follows immediately from
Proposition \ref{prop:diagH}, since \eqref{eq:cnstrts}, with $q=2s$,
is a subset of the equations that define the $\lambda_{v}$.  For the
full and restricted-full weight matrices, consider the relations in
Proposition \ref{prop:fullH} with $j \leq \tau = 2s-1$ and $i = 0$:
\begin{equation*}
  j \sum_{v=0}^{r-1} \sum_{w=0}^{r-1} h_{vw} (-1)^{j-1} (r - w)^{j-1}
  = -(-r)^{j} + \sum_{v=1}^{s} \alpha_{v} \sum_{w=0}^{v-1} \left[
  w^{j} + (w - v)^{j} \right].
\end{equation*}
Multiplying the left and right sides by $(-1)^{j-1}$, using the
symmetry of the $h_{vw}$, and swapping summation indices on the left
side, we find
\begin{equation*}
  j \sum_{v=0}^{r-1} \sigma_{v} (r - v)^{j-1} 
  = r^{j} + (-1)^{j-1}\sum_{v=1}^{s} \alpha_{v} \sum_{w=0}^{v-1} \left[
  w^{j} + (w - v)^{j} \right],
\end{equation*}
where $\sigma_{v}$ is identified with $\sum_{w=0}^{r-1} h_{vw}$.  The
second term on the right-hand side can be simplified using the
accuracy conditions of the $\alpha_{v}$ (Lemma \ref{lem:alpha}) and the
formula for the sum of powers.
\begin{align*}
  \sum_{v=1}^{s} \alpha_{v} \sum_{w=0}^{v-1} \left[
    w^{j} + (w - v)^{j} \right] 
  &= \sum_{v=1}^{s} \alpha_{v} \left[ -v^{j} 
    + \sum_{w=1}^{v} w^{j} + (-1)^j \sum_{w=1}^{v} w^{j} \right] \\
  &= \sum_{v=1}^{s} \alpha_{v} \left[ -v^{j} 
    + \frac{(1 + (-1)^{j})}{j+1} 
    \sum_{w=0}^{j} \binom{j+1}{w} \beta_{w} v^{j+1-w} \right] \\
  &= \begin{cases}
    -\frac{1}{2},         & j = 1 \\
    \phantom{-}0,         & j = 3,5,\ldots,\tau, \\
    \phantom{-}\beta_{j}, & j = 2,4,\ldots,\tau-1
    \end{cases} \\ 
  &= \beta_{j}.
\end{align*}
Thus, we have
\begin{equation*}
   j \sum_{v=0}^{r-1} \sigma_{v} (r - v)^{j-1} 
   = r^{j} - (-1)^{j} \beta_{j}, \qquad 1 \leq j \leq \tau,
\end{equation*}
and Theorem \ref{thm:eulermac} implies that full and restricted-full
SBP weight matrices are quadrature rules accurate to $\tau + 1 = 2s$.
\end{pf}

\subsection{SBP Quadrature and Coordinate Transformations}\label{sec:mapping}

Curvilinear coordinate systems are often necessary when solving PDEs
on complex domains.  Like most finite-difference schemes, SBP
operators are not applied directly to the nodes in physical space.
Instead, a coordinate transformation is used to map points in the
physical domain to points on a Cartesian grid, and the SBP operators
are applied in this uniform computational space.  However, this coordinate
transformation introduces geometric terms whose impact on the accuracy
of the quadrature rule is not clear.

We begin by considering the one-dimensional case.  Let $\fnc{T}(x) =
\xi(x)$ be an invertible transformation of class $C^{2s}$ that maps
$\Omega_{x} = [a,b]$ to $\Omega_{\xi} = [0,1]$.  For $\fnc{U} \in
L^{2}(\Omega_{x})$, the change of variable theorem implies
\begin{equation}\label{eq:int_mapping}
  \int_{a}^{b} \fnc{U}\, dx = \int_{0}^{1} \fnc{U} \fnc{J} \, d\xi, 
\end{equation}
where $\fnc{J} = \frac{dx}{d\xi}$ is the Jacobian of $\fnc{T}^{-1}$.

We are interested in the accuracy of SBP quadrature in the
computational domain, so we consider the discrete equivalent of the
right-hand side of \eqref{eq:int_mapping}.  In general the mapping will
not be explicitly available, so the Jacobian must be approximated.  As
we shall see, to retain the $2s$-order accuracy of SBP quadrature, it
is critical that the derivative that appears in the Jacobian be
approximated by the same SBP difference operator that defines the
norm.  Thus, if $x \in \mathbb{R}^{n+1}$ denotes the coordinates of the
nodes in physical space, the SBP approximation of \eqref{eq:int_mapping}
is given by
\begin{equation}\label{eq:int_sbp}
  u^{T} H D x = u^{T} Q x.
\end{equation}
The following theorem confirms that this discrete product is a
$2s$-order accurate approximation of the integral
\eqref{eq:int_mapping}.

\begin{thm}\label{thm:zdudx}
  Let $D = H^{-1} Q$ be an SBP first derivative operator.  Then
  \begin{equation*}
    (z, Du)_{H} = z^{T} Q u
  \end{equation*}
  is a $2s$-order-accurate approximation to the integral
  \begin{equation*}
    \int_{0}^{1} \fnc{Z} \frac{d \fnc{U}}{dx} \, dx,
  \end{equation*}
  where $\fnc{Z} \frac{d \fnc{U}}{dx} \in C^{2s}[0,1]$.
\end{thm}

\begin{pf}
  Using SBP-norm quadrature we have
  \begin{equation*}
    \int_{0}^{1} \fnc{Z} \frac{d \fnc{U}}{dx} \, dx = 
    (z,u')_{H} + \text{O}(h^{2s}),
  \end{equation*}
  where $u'$ denotes the analytical derivative $\partial
  \fnc{U}/\partial x$ evaluated at the grid nodes.  The result will
  follow if we can show that
  \begin{equation}\label{eq:zdudx}
    (z,u')_{H} = (z,Du)_{H} + \text{O}(h^{2s}). 
  \end{equation}
  The expression on the left of \eqref{eq:zdudx} is simply a
  quadrature for the integrand $\fnc{Y} = \fnc{Z} \frac{d}{dx} \fnc{U}$.
  Consequently, it is sufficient to show \eqref{eq:zdudx} is exact for
  polynomial integrands of degree less than $2s$.  Let 
  \begin{equation*}
    w_{i} = \begin{bmatrix} 
      x_{0}^{i} & x_{1}^{i} & \cdots & x_{n}^{i} \end{bmatrix}^{T}
  \end{equation*}
  be the restriction of the monomial $x^{i}$ to the grid.  We will
  consider
  \begin{equation*}
    z = w_{i}, \quad u = w_{j}, \quad \text{and} \qquad u' = j w_{j-1},
  \end{equation*}
  with $i+j \leq 2s$.
  
  First, suppose $j \leq s$.  In this case, an SBP operator (including
  those with diagonal-norms) is exact for $w_{j}$ giving
  \begin{equation*}
    D u = D w_{j} = j w_{j-1} = u',
  \end{equation*}
  and substitution into \eqref{eq:zdudx} yields $(z,u')_{H} =
  (z,Du)_{H}$.

  Next, to show that \eqref{eq:zdudx} is exact for $j > s$, the roles
  of $z$ and $u$ will be reversed.  Here, since $j + i \leq 2s$, we
  must have $i < s$, and the SBP operator becomes exact for $w_{i}$:
  \begin{equation*}
    D z = D w_{i} = i w_{i-1} = z'.
  \end{equation*}
  Using this exact derivative and the properties of SBP operators we find
  \begin{align*}
    (z,Du)_{H} &= z^T H (H^{-1} Q) u \\
    &= z_{n} u_{n} - z_{0} u_{0} - z^{T} Q^{T} u  \\
    &= \left.\fnc{U}\fnc{Z}\right|_{x=1} - 
       \left.\fnc{U}\fnc{Z}\right|_{x=0} - (u,Dz)_{H} \\
    &= \left.\fnc{U}\fnc{Z}\right|_{x=1} - 
       \left.\fnc{U}\fnc{Z}\right|_{x=0} - (u,z')_{H} \\
    &= \int_{0}^{1} \frac{d}{dx}\left( \fnc{U}\fnc{Z} \right)\, dx -
       \int_{0}^{1} \fnc{U} \frac{d \fnc{Z}}{dx}\, dx \\
    &= \int_{0}^{1} \fnc{Z} \frac{d \fnc{U}}{dx} \, dx.
  \end{align*}
  Thus we have shown that the expression $(z,Du)_{H}$ is also equal to
  the exact integral when $j > s$ and $i + j \leq 2s$.  This completes
  the proof.
\end{pf}

\ignore{
Applying Theorem \ref{thm:zdudx} to the integral in
\eqref{eq:int_mapping}, we have
\begin{equation*}
  \int_{a}^{b} \fnc{U}\, dx = u^{T} H (D x) + \text{O}(h^{2s}), 
\end{equation*}
which is the desired result.
}

For multidimensional problems on curvilinear tensor-product domains,
SBP operators are obtained from the one-dimensional operators using
Kronecker products.  To extend SBP quadrature to these domains, we
need only apply Theorem \ref{thm:zdudx} iteratively over the
individual coordinate directions.  We provide a sketch of the proof
here and direct the interested reader to \cite{hicken:supfun2011} for
the details of the two-dimensional case.  Consider the change of
variable theorem in $d$ dimensions:
\begin{equation*}
  \idotsint\limits_{\Omega_{x}} \fnc{W} \, dx_{1} \, dx_{2} \cdots
    dx_{d} = \idotsint\limits_{\Omega_{\xi}} \fnc{W} \fnc{J} \, d\xi_{1} \,
      d\xi_{2} \cdots d\xi_{d},
\end{equation*}
where $\fnc{J}$ is the Jacobian of the mapping (more precisely, the
determinant of the Jacobian).  As in the one-dimensional case, the
mapping and integrand must be sufficiently differentiable (class
$C^{2s}$) for the quadrature to remain $2s$-order accurate.  An
important observation is that the Jacobian consists of a sum of terms
of the form
\begin{equation*}
  \frac{\partial x_{i}}{\partial \xi_{1}} \frac{\partial x_{j}}{\partial \xi_{2}} \cdots \frac{\partial x_{k}}{\partial \xi_{d}}
\end{equation*}
in which none of the indices $i,j,\ldots,k$ are equal.  Because the 
indices of the computational coordinates are also distinct, Theorem \ref{thm:zdudx} can be applied one
dimension at a time (\ie, as an iterated integral).  For example, we
can consider dimension $\xi_{1}$ and apply Theorem \ref{thm:zdudx} to
the integral
\begin{equation*}
  \int_{0}^{1} \left( \fnc{W}  \frac{\partial x_{j}}{\partial \xi_{2}} \cdots \frac{\partial x_{k}}{\partial \xi_{d}} \right)  \frac{\partial x_{i}}{\partial \xi_{1}} \,d\xi_{1},
\end{equation*}
where $x_{i}$ corresponds with $\fnc{U}$ in the theorem, and 
\begin{equation*}
\left( \fnc{W}  \frac{\partial x_{j}}{\partial \xi_{2}} \cdots \frac{\partial x_{k}}{\partial \xi_{d}} \right) 
\end{equation*}
corresponds with $\fnc{Z}$.  Repeating this process over the remaining
coordinate directions and terms in the Jacobian yields the desired
result.

\subsection{SBP Operators and the Divergence Theorem}\label{sec:div}

Using the above results, one can show that SBP operators mimic the
$d$-dimensional divergence theorem to order $h^{2s}$ on curvilinear
domains that are diffeomorphic to the $d$-cube.  We will consider the
two-dimensional case; the extension to higher dimensions is
straightforward.

In two-dimensions, the divergence theorem is
\begin{equation}
  \iint\limits_{\Omega_{x}} \frac{\partial \fnc{F}}{\partial x} + \frac{\partial \fnc{G}}{\partial y} \, dx dy = \oint\limits_{\partial \Omega_{x}} (\fnc{F} dy + \fnc{G} dx)
\end{equation}
where $\partial \Omega_{x}$ is the piecewise-smooth boundary of
$\Omega_{x}$, oriented counter-clockwise.  Applying the coordinate
transformation, we find
\begin{align}
\iint\limits_{\Omega_{x}} \frac{\partial \fnc{F}}{\partial x} 
       + \frac{\partial \fnc{G}}{\partial y} \, dx dy 
&= \iint\limits_{\Omega_{\xi}} \left( \frac{\partial \fnc{F}}{\partial x} 
       + \frac{\partial \fnc{G}}{\partial y} \right) \fnc{J} \, dx dy \notag \\
&= \iint\limits_{\Omega_{\xi}} \frac{\partial \hat{\fnc{F}}}{\partial \xi} 
       + \frac{\partial \hat{\fnc{G}}}{\partial \eta} \, d\xi d\eta 
        \label{eq:trans_div}\\
\intertext{where we have used the metric relations \cite{pulliam:1980,thomas:1979} to obtain the components}
\hat{\fnc{F}} &= \fnc{J} \left( 
    \frac{\partial \xi}{\partial x} \fnc{F} +
    \frac{\partial \xi}{\partial y} \fnc{G} \right)
    = \phantom{-}\frac{\partial y}{\partial \eta}\fnc{F} - 
    \frac{\partial x}{\partial \eta}\fnc{G}, \label{eq:Fhat} \\
\hat{\fnc{G}} &= \fnc{J} \left( 
    \frac{\partial \eta}{\partial x} \fnc{F} +
    \frac{\partial \eta}{\partial y} \fnc{G} \right)
    = -\frac{\partial y}{\partial \xi} \fnc{F} +
    \frac{\partial x}{\partial \xi} \fnc{G}. \label{eq:Ghat}
\end{align}
In light of \eqref{eq:trans_div}, we need only show that SBP
discretizations obey the divergence theorem to order $h^{2s}$ in the
simpler computational space:
\begin{equation}\label{eq:trans_div2}
\iint\limits_{\Omega_{\xi}} \frac{\partial \hat{\fnc{F}}}{\partial \xi} 
       + \frac{\partial \hat{\fnc{G}}}{\partial \eta} \, d\xi d\eta = 
\int_{0}^{1} \left[ \hat{\fnc{F}}(1,\eta) - \hat{\fnc{F}}(0,\eta) \right]\, d\eta + 
\int_{0}^{1} \left[ \hat{\fnc{G}}(\xi,1) - \hat{\fnc{F}}(\xi,0) \right]\, d\xi. 
\end{equation}

The reader may object to this simplification, since $\hat{\fnc{F}}$
and $\hat{\fnc{G}}$ contain derivatives that depend on the geometry
and must be approximated.  However, if the partial derivatives of $x$
and $y$ appearing in \eqref{eq:Fhat} and \eqref{eq:Ghat} are
approximated using the same SBP operators as found in the discrete
divergence theorem, then Theorem \ref{thm:zdudx} can be applied.  This
follows because the same difference operator is never applied twice in the
same coordinate direction (\eg, $\partial/\partial \xi$ is applied to
$\hat{\fnc{F}}$, which contains only partial derivatives with respect
$\eta$).

For simplicity, assume that the square $\Omega_{\xi}$ is discretized
using $n+1$ nodes in both the $\xi$ and $\eta$ directions.  Thus, the
nodal coordinates are given by
\begin{equation*}
  \bsym{\xi}_{jk} = (\xi_{j},\eta_{k}) = \frac{1}{n}(j,k), 
  \qquad 0 \leq j, k \leq n \}.
\end{equation*}
If the nodes are ordered first by $j$ and then by $k$, one-dimensional
SBP operators can be used to construct the two-dimensional difference
operators
\begin{gather*}
  D_{\xi} = (I \otimes D), \qquad \text{and}\quad 
  D_{\eta} = (D \otimes I),
\end{gather*}
where $\otimes$ denotes the Kronecker product, $D = H^{-1} Q$ is the
one-dimensional SBP operator, and $I$ is the $(n+1)\times (n+1)$
identity matrix.  Similarly, $(H \otimes H)$ defines the SBP
quadrature for the two-dimensional set of points.  Let $B =
\diag{(-1,0,0,\ldots,1)}$, so that we may write $Q + Q^{T} = B$.
Finally, let $\hat{f}$ and $\hat{g}$ denote the restriction of the
functions $\hat{\fnc{F}}$ and $\hat{\fnc{G}}$, respectively, to the
grid points, and let $c = \begin{bmatrix} 1 & 1 & \cdots &
  1 \end{bmatrix}^{T}$ denote the constant function $1$ restricted to
the grid.

With the two-dimensional SBP operators suitably defined, we can
discretize the left-hand side of \eqref{eq:trans_div2}:
\begin{align}
  c^{T} (H \otimes H) \left[ (I \otimes D) \hat{f} + (D \otimes I) \hat{g} \right] &= c^{T} (H \otimes Q) \hat{f} + c^{T} (Q \otimes H) \hat{g} \notag \\
   &= c^{T} \left( H \otimes (B - Q^{T}) \right) \hat{f} 
  + c^{T} \left((B - Q^{T}) \otimes H \right) \hat{g} \notag \\
   &= \sum_{j=0}^{n} h_{ii} (\hat{f}_{n,j} - \hat{f}_{0,j}) + \sum_{i=0}^{n} h_{ii} (\hat{g}_{i,n} - \hat{g}_{i,0}), \label{eq:sbp_div}
\end{align}
where we have used $c^{T} (Q^{T} \otimes H) = c^{T} (H \otimes Q^{T})
= 0$ (constants are in the null space of $D = H^{T} Q$).  

We highlight two significant facts regarding \eqref{eq:sbp_div}.
\begin{enumerate}
\item It is a $2s$-order accurate approximation of the right-hand
side of \eqref{eq:trans_div2}.  
\item It depends only on the terms of $\hat{f}$ and
  $\hat{g}$ that fall on the boundary.
\end{enumerate}
Constructing a scheme that satisfies either one of these properties
may not be difficult; however, few high-order schemes satisfy both 1
and 2 simultaneously.  This is what we mean when we say the SBP
operator mimics the divergence theorem.

\section{Examples}\label{sec:example}

\subsection{One-dimensional Quadrature}

To illustrate the basic theory, we use the weight matrices from several
common SBP operators to integrate a simple function.  We consider
three SBP operators with diagonal weight matrices and one SBP operator
with a full norm.  The diagonal operators are taken from Diener \etal
\cite{diener:2007} and are denoted by diag-$\tau$-$2s$, where $\tau$
and $2s$ indicate the truncation error at the boundary and interior,
respectively.  The full norm operator can be found in
\cite{strand:1994} and is denoted full-$\tau$-$2s$.  The boundary
weights $\sigma_{v}$ for all four operators are listed in Table
\ref{tab:SBPrules}; for the diagonal norms $\sigma_{v} = \lambda_{v}$,
whereas for the full norm $\sigma_{v} = \sum_{w=0}^{r-1}h_{vw}$.

\begin{table}[tbp]
  \begin{center}
    \caption[]{\small Boundary quadrature weights corresponding to
    some SBP weight matrices.}\label{tab:SBPrules}
    \begin{threeparttable}
      \begin{tabular}{ccccccccc}\hline
	\rule{0ex}{3ex}\textbf{SBP operator} & $\bsym{\tau}$ & $\mathbf{2s}$
	& $\bsym{\sigma}_{0}$ & $\bsym{\sigma}_{1}$ 
	& $\bsym{\sigma}_{2}$ & $\bsym{\sigma}_{3}$ 
	& $\bsym{\sigma}_{4}$ & $\bsym{\sigma}_{5}$ \\\hline 
	\rule{0ex}{3ex}diag-1-2\tnote{\ensuremath{\dagger}} & 1 & 2 & 
	$\frac{1}{2}$ & --- & --- & --- & --- & --- \\
	\rule{0ex}{3ex}diag-2-4 & 2 & 4 &
	$\frac{17}{48}$ & $\frac{59}{48}$ & $\frac{43}{48}$ & $\frac{49}{48}$ &
	--- & --- \\
	\rule{0ex}{3ex}full-3-4 & 3 & 4 &
	$\frac{43}{144}$ & $\frac{67}{48}$ & $\frac{35}{48}$ &
	$\frac{155}{144}$ & --- & --- \\
	\rule{0ex}{3ex}diag-3-6 & 3 & 6 &
	$\frac{13649}{43200}$ & $\frac{12013}{8640}$ & $\frac{2711}{4320}$ &
	$\frac{5359}{4320}$   & $\frac{7877}{8640}$  & $\frac{43801}{43200}$
	\\[1ex]\hline
      \end{tabular}
      \begin{tablenotes}
      \item[\ensuremath{\dagger}] {\small the trapezoidal rule}
      \end{tablenotes}
    \end{threeparttable}
  \end{center}
\end{table}

Consider the definite integral
\begin{align}
   \fnc{I} 
   &= \int_{0}^{1} \fnc{U}(x) \, dx \notag\\ 
   &= \int_{0}^{1} (4 \pi)^{2} x \cos{(4 \pi x)} \, dx \label{eq:testint} \\
   &= - 4\pi \cos{(4\pi)}. \notag
\end{align}
To assess the accuracy of the SBP quadrature rules in
Table~\ref{tab:SBPrules}, we perform a grid refinement study based on
the integral~\eqref{eq:testint} and using $n \in \{16, 32, 64, 128,
256, 512 \}$.  Table~\ref{tab:quaderror} lists the rates of
convergence for the quadrature rules.  For $n > 16$, the rate of
convergence is calculated from
\begin{equation}\label{eq:rate}
  q_{n} = \frac{1}{\ln{(2)}} \ln{ \left( \frac{|E_{\frac{n}{2}}|}{|E_{n}|}\right)},
\end{equation}
where $E_{n} = \fnc{I} - c^{T} H u$, with $c^{T} \equiv \begin{pmatrix} 1 & 1 & \ldots & 1 \end{pmatrix}$, is the error using $n+1$
nodes.  In all cases, the errors converge to zero at the expected
asymptotic rate of $2s$.

\begin{figure}[tbp]
  \begin{center}
    \includegraphics[width=0.8\textwidth]{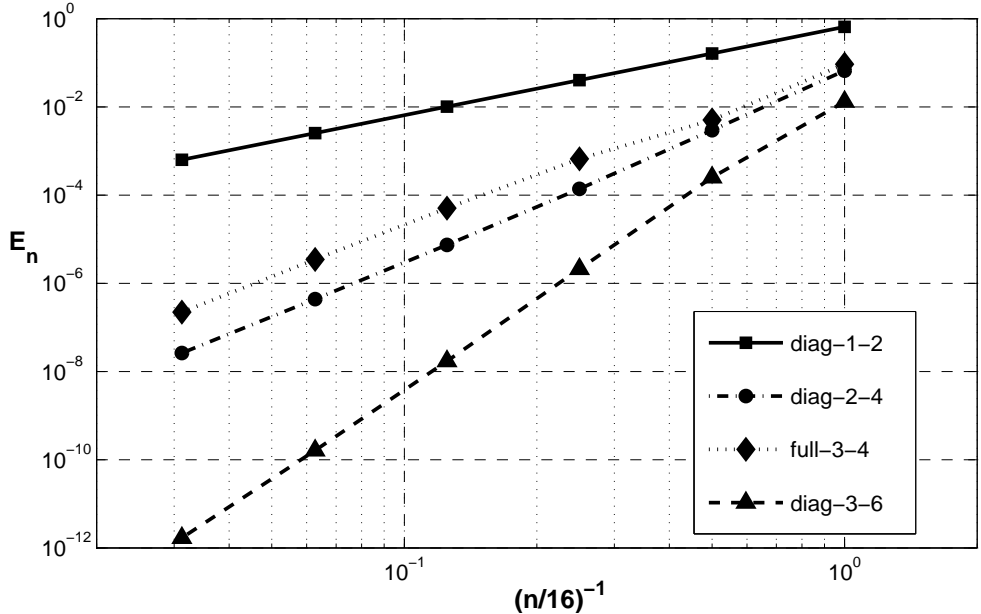}
    \caption[]{ Errors of the SBP-based quadrature rules applied to
    \eqref{eq:testint}.\label{fig:quaderrors}}
  \end{center}
\end{figure}

Figure~\ref{fig:quaderrors} plots the errors $E_{n}$ versus a
normalized mesh spacing.  This figure reminds us that schemes with the
same order of accuracy can produce different absolute errors: the
diag-2-4 operator is almost an order of magnitude more accurate than
the full-3-4 operator for $n \geq 64$.  However, further analysis is
required before we can characterize the relative performance of these
schemes more generally.

\begin{table}[tbp]
  \begin{center}
    \caption[]{\small Rates of convergence for the SBP quadrature
    rules in Table \ref{tab:SBPrules} applied to
    \eqref{eq:testint}.}\label{tab:quaderror}
    \begin{threeparttable}
      \begin{tabular}{cccccc}\hline
	 & \multicolumn{5}{c}{$\mathbf{n}$} \\\cline{2-6}
	\textbf{SBP operator}\rule{0ex}{3ex}
	& 32 & 64 & 128 & 256 & 512 \\\hline
	\rule{0ex}{3ex}diag-1-2 &
	2.0113 & 2.0028 & 2.0007 & 2.0002 & 2.0000 \\
	\rule{0ex}{3ex}diag-2-4 &
	4.4978 & 4.4148 & 4.2182 & 4.1019 & 4.0473 \\
	\rule{0ex}{3ex}full-3-4 &
        4.1973 & 2.9369 & 3.7072 & 3.8876 & 3.9510 \\
	\rule{0ex}{3ex}diag-3-6 &
	5.7050 & 6.8942 & 6.9378 & 6.7651 & 6.5472 \\\hline
      \end{tabular}
      %\begin{tablenotes}
      %\item[\ensuremath{\dagger}] {\small also an SBP diagonal norm}
      %\end{tablenotes}
    \end{threeparttable}
  \end{center}
\end{table}

\subsection{Multi-dimensional Quadrature on a Curvilinear Domain}

\begin{figure}[tbp]
  \begin{center}
    \includegraphics[width=0.5\textwidth]{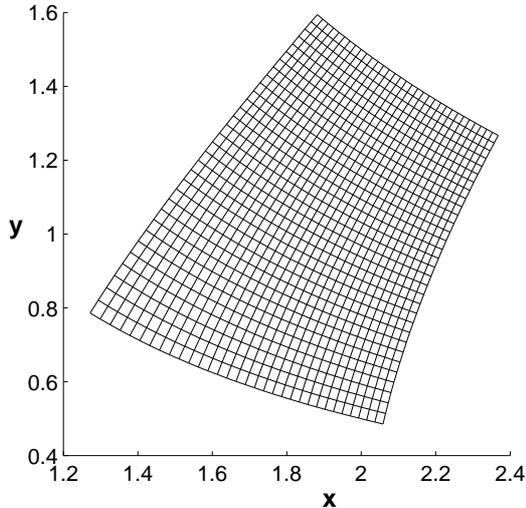}
    \caption[]{Example grid for $\Omega_{x}$ with $n = 32$.
      \label{fig:mapdomain}}
  \end{center}
\end{figure}

As shown in Section \ref{sec:mapping}, SBP quadrature retains its
theoretical accuracy on curvilinear domains provided the Jacobian of
the transformation is approximated using the corresponding SBP
difference operator.  To verify this, we consider the domain
\begin{align}
  \Omega_{x} &= \{ (x,y) \in \mathbb{R}^{2} \mid 1 \leq xy \leq 3, 1 \leq x^{2} - y^{2} \leq 4 \}, \notag \\
\intertext{and the integral}
  \fnc{I} &= \iint\limits_{\Omega_{x}} (x^{2} + y^{2}) e^{\frac{1 - x^{2} + y^{2}}{3}} \sin{\left(\frac{xy - 1}{2}\right)} \, dx dy \notag \\
          &= 3(1 - e^{-1})(1 - \cos(1)). \label{eq:mapint}
\end{align}
To compute this integral numerically, we introduce a computational
domain based on the coordinates
\begin{gather*}
  \xi = \frac{x^{2} - y^{2} - 1}{3}, \qquad \text{and} \quad \eta = \frac{xy -1}{2}.
\end{gather*}
For a given $n \in \{16, 32, 64, 128, 256, 512 \}$, we divide $\xi$
and $\eta$ uniformly into $n+1$ points to produce a Cartesian grid on
the square $\Omega_{\xi} = [0,1]^{2}$.  The physical coordinates $x$
and $y$ are evaluated at each computational coordinate, and these are
used to compute the integrand in \eqref{eq:mapint}, which we denote by
$f$.  The grid for $n = 32$ is shown in Figure \ref{fig:mapdomain}.

The Jacobian of the transformation is approximated using
\begin{equation}\label{eq:Jac2d}
  J = \left[ (I \otimes D) x \right] \circ \left[ (D \otimes I) y \right] 
  - \left[ (I \otimes D) y \right] \circ \left[ (D \otimes I) x \right], 
\end{equation}
where $\circ$ denotes the Hadamard product (the entry-wise product,
analogous to matrix addition).  We have assumed that the nodes are
ordered first by $\xi$ and then by $\eta$, so we can construct the
two-dimensional derivative operators using Kronecker products of the
one-dimensional operator $D$ and identity matrix $I$.

\begin{figure}[tbp]
  \begin{center}
    \includegraphics[width=0.8\textwidth]{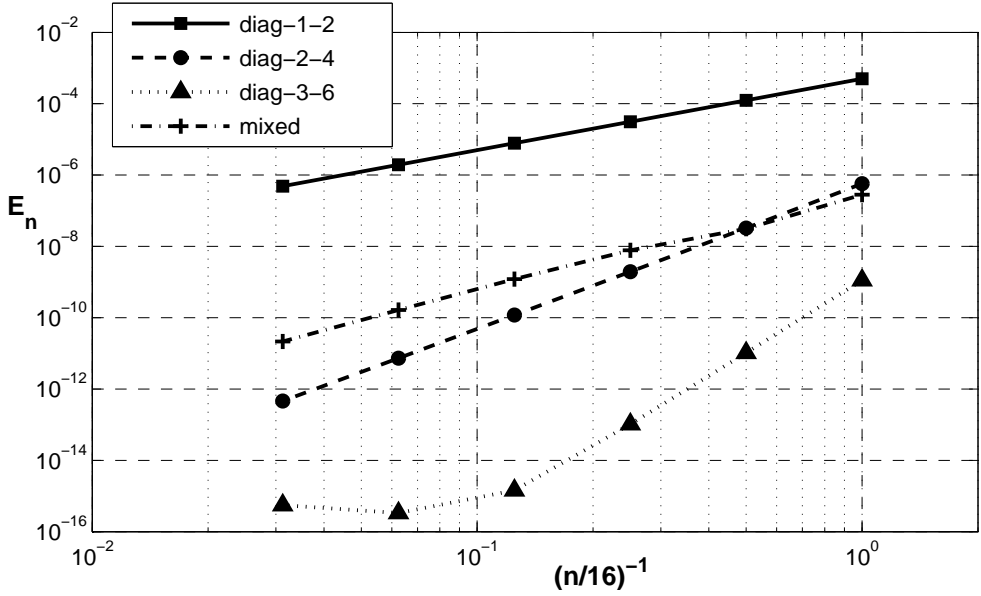}
    \caption[]{Errors of the SBP-based quadrature rules in
      approximating \eqref{eq:mapint}.\label{fig:mapping}}
  \end{center}
\end{figure}

For a given $n$, the SBP-based approximation of \eqref{eq:mapint} is given by
\begin{equation} \label{eq:mapsbp}
  \fnc{I}_{n} = J^{T} (H \otimes H) f, 
\end{equation}
and the error in the quadrature is $E_{n} = \fnc{I} - \fnc{I}_{n}$.
As before, the order of convergence for $n > 16$ is estimated by
$q_{n}$ given by \eqref{eq:rate}.  Figure \ref{fig:mapping} plots
$E_{n}$ and Table \ref{tab:mapping} lists $q_{n}$ for the
diagonal-norm SBP operators listed in Table \ref{tab:SBPrules}.  We
focus on the diagonal-norm SBP operators here, because their accuracy
is less obvious; the derivatives used to approximate the Jacobian are
only $s$-order accurate at the boundary.  Nevertheless, as predicted
by the theory, Table \ref{tab:mapping} shows that the quadrature for
these schemes remains $2s$-order accurate.  Note that the errors for
the diag-3-6 scheme are corrupted by round-off error for $n = 256$ and
$n = 512$, which explains the suboptimal values of $q_{n}$ for these
grids.

We have also included results for a mixed scheme in Table
\ref{tab:mapping} and Figure \ref{fig:mapping}.  This mixed scheme
uses the diag-3-6 SBP operator to evaluate the derivatives in the
Jacobian \eqref{eq:Jac2d} and the diag-2-4 operator to evaluate the
quadrature \eqref{eq:mapsbp}.  The results show that the mixed scheme
has an asymptotic convergence rate of only 3.  Thus, despite a more
accurate approximation of the Jacobian, the mixed scheme produces a
less accurate $\fnc{I}_{n}$ than the scheme using the diag-2-4
operator for both the Jacobian and quadrature.  This illustrates the
importance of using the same operator to obtain the theoretical
convergence rate.

\begin{table}[tbp]
  \begin{center}
    \caption[]{\small Rates of convergence for the diagonal-norm SBP operator
    approximation of \eqref{eq:mapint}.}\label{tab:mapping}
    \begin{threeparttable}
      \begin{tabular}{cccccc}\hline
	 & \multicolumn{5}{c}{$\mathbf{n}$} \\\cline{2-6}
	\textbf{SBP operator}\rule{0ex}{3ex}
	& 32 & 64 & 128 & 256 & 512 \\\hline
	\rule{0ex}{3ex}diag-1-2 &
	 2.0911 & 2.0453 & 2.0226 & 2.0113 & \phantom{-}2.0056 \\
	\rule{0ex}{3ex}diag-2-4 &
	 4.3283 & 4.1583 & 4.0768 & 4.0374 & \phantom{-}4.0093 \\
	\rule{0ex}{3ex}diag-3-6 &
	 7.0799 & 6.7941 & 6.2253 & 2.1274 & -0.7390 \\
         \rule{0ex}{3ex}mixed &
         3.3170 & 2.0521 & 2.7215 & 2.8863 & \phantom{-} 2.9484 \\\hline
      \end{tabular}
      %\begin{tablenotes}
      %\item[\ensuremath{\dagger}] {\small also an SBP diagonal norm}
      %\end{tablenotes}
    \end{threeparttable}
  \end{center}
\end{table}

\subsection{Discrete Divergence Theorem}

In this final example, we verify that SBP operators mimic the
divergence theorem accurately.  Specifically, we wish to show that
when the divergence of a vector field is discretized using SBP
operators and then integrated using the corresponding SBP quadrature
rule, the result depends only on the nodes along the boundary and is a
$2s$-order approximation to the surface flux.

We adopt the same domain $\Omega_{x}$ and coordinate transformation
as in the previous example.  A vector field $(\fnc{F},\fnc{G})$ is defined by 
\begin{align*}
  \fnc{F}(x,y) &= \phantom{-}\frac{x}{2} \exp{\left(\frac{1 - xy}{2}\right)} \cos{\left( \frac{2 \pi(x^{2} - y^{2} - 1)}{3} \right)} + \frac{2y}{3} \left( \frac{xy-1}{2} \right)^{7} \sin{\left( \frac{\pi(x^{2} - y^{2} - 1)}{3} \right)} \\
  \fnc{G}(x,y) &= -\frac{y}{2} \exp{\left(\frac{1 - xy}{2}\right)} \cos{\left( \frac{2 \pi(x^{2} - y^{2} - 1)}{3} \right)} + \frac{2x}{3} \left( \frac{xy-1}{2} \right)^{7} \sin{\left( \frac{\pi(x^{2} - y^{2} - 1)}{3} \right)}.
\end{align*}
The analytical value of the divergence of $(\fnc{F},\fnc{G})$
integrated over the domain $\Omega_{x}$ is 
\begin{equation}
  \fnc{I} = \iint\limits_{\Omega_{x}} \frac{\partial \fnc{F}}{\partial x} + \frac{\partial \fnc{G}}{\partial y} \, dx dy = \frac{2}{\pi}. \label{eq:mapdiv}
\end{equation}

The discrete divergence is evaluated in computational space using
approximations for $\hat{\fnc{F}}$ and $\hat{\fnc{G}}$.  In
particular, the derivatives of the spatial coordinates that appear in
\eqref{eq:Fhat} and \eqref{eq:Ghat} are approximated using SBP
operators.  Therefore, at the nodes, $\hat{\fnc{F}}$ and
$\hat{\fnc{G}}$ take on the values
\begin{align*}
  \hat{f} &= \phantom{-}\left[ (D \otimes I) y \right] \circ f 
             - \left[ (D \otimes I) x \right] \circ g, \\
  \hat{g} &= -\left[ (I \otimes D) y \right] \circ f 
             + \left[ (I \otimes D) x \right] \circ g, \\
\end{align*}
where $f$ and $g$ denote the values of $\fnc{F}$ and $\fnc{G}$
evaluated at the nodes.

The SBP approximation of $\fnc{I}$ is given by (see \eqref{eq:trans_div})
\begin{equation*}
  \fnc{I}_{n} = c^{T} (H \otimes H) \left[ (I \otimes D) \hat{f} + (D \otimes I) \hat{g} \right].
\end{equation*}
Table \ref{tab:div} lists the estimated order of accuracy $q_{n}$
based on $E_{n} = \fnc{I} - \fnc{I}_{n}$ for the three diagonal-norm
SBP operators diag-1-2, diag-2-4, and diag-3-6.  As predicted, the
SBP discrete divergence integrated using $(H \otimes H)$ is a
$2s$-order accurate approximation to $\fnc{I}$.  Moreover, in light of
\eqref{eq:sbp_div}, we know that $\fnc{I}_{n}$ depends only on the
boundary nodes (This has been confirmed by calculating the right-hand side
of \eqref{eq:sbp_div} and showing that it equals $\fnc{I}_{n}$ to
machine error).

\begin{table}[tbp]
  \begin{center}
    \caption[]{\small Rates of convergence for the diagonal-norm SBP
      operator approximation of an integrated divergence field.
      Round-off errors are contaminating the estimates for diag-3-6
      with $n = 256$ and $n = 512$}\label{tab:div}
    \begin{threeparttable}
      \begin{tabular}{cccccc}\hline
	 & \multicolumn{5}{c}{$\mathbf{n}$} \\\cline{2-6}
	\textbf{SBP operator}\rule{0ex}{3ex}
	& 32 & 64 & 128 & 256 & 512 \\\hline
	\rule{0ex}{3ex}diag-1-2 &
	2.0909 & 2.0453 & 2.0226 & 2.0113 & \phantom{-}2.0056 \\
	\rule{0ex}{3ex}diag-2-4 &
	3.7201 & 3.7862 & 3.9000 & 3.9532 & \phantom{-}3.9758 \\
	\rule{0ex}{3ex}diag-3-6 &
	7.5935 & 7.2371 & 7.8361 & 5.0507 & -2.1760
 \\\hline
      \end{tabular}
      %\begin{tablenotes}
      %\item[\ensuremath{\dagger}] {\small also an SBP diagonal norm}
      %\end{tablenotes}
    \end{threeparttable}
  \end{center}
\end{table}

\section{Discussion and Conclusions}\label{sec:conclude}

We have shown that the weight matrices of SBP finite-difference
operators are related to trapezoid rules with end corrections.  We
make no claim regarding the optimality of these SBP quadrature rules
with respect to existing schemes on uniformly spaced grid points.
However, the result has significant implications for SBP
discretizations of PDEs, which we list below.
\begin{itemize}
\item The SBP energy norm, which is frequently used in the stability
  analysis of SBP-based PDE discretizations, is a $\text{O}(h^{2s})$
  accurate approximation of the $L^{2}$ norm for functions on $[0,1]$.
\item The summation-by-parts property, equation \eqref{eq:SBP}, is a
  formal \emph{and} accurate representation of integration by parts,
  equation \eqref{eq:IBP}.  More generally, multi-dimensional SBP
  discretizations using Kronecker products mimic the divergence
  theorem, \ie the weight-matrix quadrature applied to the discrete
  divergence produces an accurate quadrature of the flux over the
  domain boundary in which no interior points are involved.
\item Diagonal-norm SBP operators have $s$ order-accurate boundary
  closures when the interior scheme is $2s$-order accurate.  This
  limits numerical PDE solutions to $s+1$ order accuracy
  \cite{gustafsson:1975}; however, in \cite{hicken:supfun2011} we show
  that a dual consistent SBP discretization leads to super-convergent
  $2s$-order-accurate functionals, if the corresponding SBP quadrature
  rule is used to calculate the functional (see also
  \cite{hicken:cfdsc2010}).
\end{itemize}
In light of these observations, the SBP weight matrix appears to be the
natural quadrature rule for evaluating functionals from corresponding
SBP discretizations.

Other than its asymptotic form, the leading error term for SBP-based
quadrature remains unknown, so we cannot make general statements
regarding the relative accuracy of different weight matrices.
Although not pursued here, characterizing the leading error term
should be possible by finding the finite-difference schemes
corresponding to the $\delta_{v}^{(k-1)}$ in Theorem
\ref{thm:eulermac}.

\ignore{
We conclude with a conjecture.  The SBP first derivative definition
involves three matrices: the difference operator, $D$, the weight
matrix $H$, and the matrix $Q$.  The significance of $Q$ has not been
discussed, but we postulate that it represents a one-dimensional
finite-volume operator, where the cell volumes are defined by the
$\sigma_{v}$ in $H$.  A thorough analysis of $Q$ will be the subject
of a future paper.}

\bibliographystyle{elsarticle-num}
\bibliography{/home/jehicken/Biblio/ref} 

\end{document}